\documentclass[11pt,twoside]{amsart}
\newtheorem{thm}{Theorem}[section]

\newtheorem{lem}[thm]{Lemma}

\newtheorem{rem}[thm]{\bf{Remark}}

\numberwithin{equation}{section}
\def\pn{\par\noindent}

\begin{document}

\leftline{ \scriptsize}

\vspace{1.3 cm}

\title
{A note on a paper of Harris concerning the asymptotic approximation to the eigenvalues of $-y''+qy=\lambda y$, with boundary conditions of general form}
\author{Mahdi Hormozi}

\thanks{{\scriptsize
\hskip -0.4 true cm MSC(2000): Primary: 41A05, 34B05 ; Secondary:94A20
$*$Corresponding author}}
\maketitle

\begin{center}

\end{center}

\begin{abstract}  In this paper, we derive an asymptotic approximation to the eigenvalues of the linear differential equation
$$
-y''(x)+q(x)y(x)=\lambda y(x),\hskip 1.4 true cm x\in (a,b)
$$
with boundary conditions of general form, when $q$ is a measurable function which has a singularity in $(a,b)$ and which is integrable on subsets of $(a,b)$ which exclude the singularity.
\end{abstract}

\vskip 0.2 true cm


\pagestyle{myheadings}
\markboth{\rightline {\scriptsize  Mahdi Hormozi}}
         {\leftline{\scriptsize }}

\bigskip
\bigskip


\vskip 0.4 true cm

\section{\bf Introduction}

\vskip 0.4 true cm
Consider the linear differential equation
\begin{equation}
-y''(x)+q(x)y(x)=\lambda y(x),\hspace{0.5cm} x\in (a,b),
\end{equation}
where $\lambda$ is a real parameter and $q$ is real-valued function which has a singularity in $(a,b)$. According to [2], an eigenvalue problem may be associate with (1.1) by imposing the boundary conditions
 \begin{eqnarray}
y(a) \cos\alpha  - y'(a)\sin \alpha =0 ,\hspace{0.5cm} \alpha \in [0,\pi),\\
 y(b) \cos \beta   - y'(b)\sin \beta =0,\hspace{0.5cm} \beta \in [0,\pi).
 \end{eqnarray}

 In [1], Atkinson obtained an asymptotic approximation of eigenvalues where $y$ satisfies Dirichlet and Neumann boundary conditions in (1.1). Here, we find asymptotic approximation of eigenvalues for all boundary condition of the forms (1.2) and (1.3). To achieve this, we transform (1.1) to a differential equation all of whose coefficients belong to $L_{1}[a, b]$. Then
we employ a Pr\"{u}fer transformation to obtain an approximation of the eigenvalues. In this way, many basic properties of singular problems can be inferred from the corresponding regular
ones. In [8], Harris derived an asymptotic approximation to the eigenvalues of the differential equation (1.1), defined on the interval $[a,b]$, with boundary conditions of general form. But, he demands the condition, $q\in L^1[a,b]$. Atkinson and Harris found asymptotic formulae for the eigenvalues of spectral problems associated with linear differential equations of the form (1.1), where $q(x)$ has a singularity of the form $\alpha x^{-k}$ with $1\leq k<\tfrac43$ and $1\leq k<\tfrac32$ in [1] and [7] respectively. Harris and Race generalized those results for the case $1\leq k<2$ in [10]. In [11], Harris and  Marzano derived asymptotic estimates for the eigenvalues of (1.1) on $[0,a]$ with periodic and semi-periodic boundary conditions.
 We consider $q(x)= Cx^{-K}$ where $ 1 \leq K<2$ and an asymptotic approximation to the eigenvalues of (1.1) with boundary conditions of general form. Our technique in this paper follows closely the  technique used in [1],[7-8] and [10]. Let $U= [a,0)\cup (0,b]$ and $q \in L_{1,Loc}(U)$. As Harris did in [10, pp. 90], suppose that there exists some real function $f$ on $[a,0)\cup (0,b]$ in $AC_{Loc}([a,0)\cup (0,b])$ which regularizes (1.1) in the following sense. For $f$ which can be chosen in section 2, define quasi-derivatives, $y^{[i]}$ as follows:

$$y^{[0]}:=y,\hspace{0.5cm} y^{[1]}:=y'+fy,$$
$y$ is a solution of (1.1) with boundary conditions (1.2) and (1.3) if and only if
\begin{equation}
\left( {\begin{array}{*{20}c}
   y^{[0]} \\
   y^{[1]} \\

 \end{array} } \right)= \left( {\begin{array}{*{20}c}
   -f & 1  \\
   f'+q-f^{2}- \lambda & f  \\

 \end{array} } \right)\left( {\begin{array}{*{20}c}
   y^{[0]} \\
   y^{[1]} \\

 \end{array} } \right)
  \end{equation}
The object of the regularization process is to chose $f$ in such way that
\begin{equation}
 f \in L^1(a,b) \hspace{0.5cm} \textmd{and }\hspace{0.5cm} -F:=q-f^2+f' \in  L^1(a,b).
\end{equation}
Having rewritten (1.1) as the system (1.4), we observe that, for any solution $y$ of (1.1) with
$\lambda > 0$, according to [1] and [7], we can define a function $ \theta \in AC(a, b)$ by
$$
\tan \theta = \frac{\lambda ^{\frac{1}{2}} y}{y^{[1]}}.
$$
When $y^{[1]}=0$, $\theta$ is defined by continuity [10, p. 91]. It makes sense to mention that one can find full discussions and nice examples about the choice of $f$ in [1], [7] and [10]. Atkinson in [1] noticed that the function $\theta$ satisfies the differential equation

\begin{equation}
\theta'= \lambda ^{\frac{1}{2}}-f\sin(2\theta)+ \lambda ^{-\frac{1}{2}}F\sin^2(\theta).
\end{equation}
Let $\lambda >0$ and the $n$-th eigenvalue $\lambda_n$ of (1.1-3), then according to [2, theorem 2], Dirichlet and non-Dirichlet boundary conditions can be described as bellow:
$$
\
\begin{cases}
\text{in Case 1}~(\alpha=0,\beta =0):& \theta (b,\lambda)-\theta (a,\lambda)=(n+1)\pi;\\[.1in]
\text{in Case 2}~(\alpha=0,\beta \neq 0):& \theta (b,\lambda)-\theta (a,\lambda)=(n+\frac{1}{2})\pi-\lambda^{-\frac{1}{2}}\cot \beta + O(\lambda^{-\frac{3}{2}});\\[.1in]
\text{in Case 3}~ (\alpha\neq 0,\beta =0):& \theta (b,\lambda)-\theta (a,\lambda)=(n+\frac{1}{2})\pi+\lambda^{-\frac{1}{2}}\cot \alpha + O(\lambda^{-\frac{3}{2}});\\[.1in]
\text{in Case 4}~ (\alpha\neq0,\beta \neq0):& \theta (b,\lambda)-\theta (a,\lambda)=n\pi+\lambda^{-\frac{1}{2}}(\cot \alpha-\cot \beta) + O(\lambda^{-\frac{3}{2}}).\\[.1in]
\end{cases}
$$
It follows from (1.5-6) that large positive eigenvalues of either the Dirichlet or non-Dirichlet problems over $[a,b]$ satisfy
\begin{equation}
\lambda ^{\frac{1}{2}}=  \frac{\theta(b)-\theta(a)}{(b-a)}+ O(1).
\end{equation}
Our aim here is to obtain a formula like (1.7) in which the $O(1)$ term is replaced by an integral term plus and error term of smaller order. We obtain an error term of $o(\lambda^{\frac{-N}{2}}) (N\geq1)$. To achieve this we first use the differential equation (1.6) to obtain estimates for $\theta(b)-\theta(a)$ for general $\lambda$ as $\lambda \rightarrow \infty$.
\vskip 0.4 true cm

\section{Statement of Result}
\vskip 0.4 true cm
We define a sequence ${\xi_{j}(t)}$ for $j=1,...,N+1,~t \in [a,b]$ by\\
\begin{eqnarray}
 \nonumber
\xi_{1}(t) := |\int_{0}^{t} |f(s)| + |F(s)| ds|\\
\\
 \nonumber
\xi_{j}(t) := |\int_{0}^{t}(|f(s)| + |F(s)|)\xi_{j-1}(s) ds|
\end{eqnarray}
and note that in view of $f,F\in L(a,b)$ ,
\begin{eqnarray}
\xi_{j}(t)\leq c \xi_{j-1}(t)\hspace{0.5cm}\textmd{for}  \hspace{0.5cm} t \in [a,b], \hspace{0.5cm}  2 \leq j\leq N+1
\end{eqnarray}

Suppose that for some $N \geq 1$,\\
\begin{eqnarray}
 \nonumber
f'\xi_{N+1},\hspace{0.5cm} f^{2}\xi_{N}, \hspace{0.5cm}fF\xi_{N} \in L[a,b];\\
\\
 \nonumber
f(t)\xi_{N+1}(t) \rightarrow 0  \hspace{0.5cm} \textmd{as} \hspace{0.5cm} t \rightarrow 0.
\end{eqnarray}

\vskip 0.4 true cm

We define a sequence of approximating functions a

\begin{equation}
\theta_0(x):=\theta(a)+ \lambda ^{\frac{1}{2}}(x-a);
\end{equation}

\begin{equation}
 \theta_{j}(0):=\theta(0);
\end{equation}

\begin{equation}
\theta_{j+1}(x):=\theta(a)+ \lambda ^{\frac{1}{2}}(x-a)-\int_{a}^{x}f\sin(2\theta_{j}(t))dt + \lambda ^{-\frac{1}{2}}\int_{a}^{x}F\sin^2(\theta_{j}(t))dt.
\end{equation}

for $j=0,1,2,...$ and for $a\leq x \leq b.$ We measure the closeness of the approximation in the next result. Thus
\begin{equation}
\theta'_{j+1}= \lambda ^{\frac{1}{2}}-f\sin(2\theta_{j})+ \lambda ^{-\frac{1}{2}}F\sin^2(\theta_{j})
\end{equation}
The following lemma appears in [1] and [10].
\vskip 0.4 true cm
\begin{lem}If $g\in \L^1$ then for any $j$ and $ a \leq x\leq b$
$$
\int_{a}^{x}g(t) \sin(2\theta_j(t))dt =o(1)
$$
 as $\lambda \rightarrow \infty$.
\end{lem}

\vskip 0.4 true cm
By using lemma 5.1 and lemma 5.2 of [10] we conclude the following lemma
\begin{lem} There exists a suitable constant $C$ such that
$$| \theta_{j+1}-\theta_{j}|\leq C \sup_{ a\leq x\leq b} | \theta-\theta_{j}|\xi_{j+1}(x)$$
\end{lem}

\vskip 0.4 true cm
 Now, we prove an elementary lemma.
\vskip 0.4 true cm
\begin{lem}If $g\in \L^1$ and $\theta(x)- \theta_{j}(x)= \lambda ^{-\frac{1}{2}} \int_{a}^{x}g\{\sin^2(\theta(t)) - \sin^2(\theta_{j}(t)) \}dt $
 then $|\theta(x)-\theta_{j+1}(x) |\leq \lambda ^{-\frac{1}{2}} \sup_{ a\leq x\leq b}|\theta(x)-\theta_{j}(x) |\int_{a}^{x}gdt  $
 \end{lem}
\pn{\bf Proof.}
\begin{eqnarray*}
\theta(x)-\theta_{j+1}(x) &=&\lambda ^{-\frac{1}{2}} \int_{a}^{x}g\{\sin^2(\theta(t)) - \sin^2(\theta_{j}(t))           \}dt \\
&=& \frac{1}{2}\lambda ^{-\frac{1}{2}} \int_{a}^{x}g\{\cos(2\theta_j(t)) - \cos(2\theta(t))           \}dt\\
&=&-\lambda ^{-\frac{1}{2}} \int_{a}^{x}g \sin(\theta_j(t)-\theta(t))\sin(\theta_j(t)+\theta(t))           \}dt\\
&\leq&\lambda ^{-\frac{1}{2}} \sup_{ a\leq x\leq b}|\theta(x)-\theta_{j}(x) |\int_{a}^{x}gdt
\end{eqnarray*}
\vskip 0.4 true cm
\begin{rem} Lemma 2.2 shows that if $|\theta(x)- \theta_{j}(x)|= o(\lambda^{\frac{-j}{2}})$ then $|\theta(x)- \theta_{j+1}(x)|= o(\lambda^{\frac{-(j+1)}{2}})$
\end{rem}
\vskip 0.4 true cm
\begin{lem}There exists a suitable constant $C$ such that
$$\int_{a}^{x}f \left(\sin(2\theta_{j}(t))- \sin(2\theta(t))\right)dt \leq C\lambda ^{-\frac{1}{2}} \sup_{ a\leq x\leq b}|\theta(x)-\theta_{j}(x)|\hspace{0.5cm}, x\in (a,b),$$
\end{lem}
\vskip 0.4 true cm
\pn{\bf Proof.}

\begin{eqnarray*}
\int_{a}^{x}f \left(\sin(2\theta(t))(1) - \sin(2\theta_j(t))(1)\right) dt &=&\lambda ^{-\frac{1}{2}} \int_{a}^{x} f \{ \sin(2\theta)\theta'  -    \sin(2\theta_{j}) {\theta_j}'     \} dt\\
&+& \lambda ^{-\frac{1}{2}} \int_{a}^{x} f^2 \{ \sin^2(2\theta) -  \sin(2\theta_j) \sin(2\theta_{j-1})     \} dt\\
&-& \lambda ^{-1} \int_{a}^{x} fF \{ \sin(2\theta)\sin^2(\theta) - \sin(2\theta_j)\sin^2(\theta_{j-1})      \} dt\\
&=:&I_1 + I_2 - I_3.
\end{eqnarray*}
But
$$I_1= \lambda ^{\frac{-1}{2}}[f(t)( \sin^2(\theta(t)) - \sin^2(\theta_j (t)) ) ]_{a}^{x}-\lambda ^{-\frac{1}{2}} \int_{a}^{x} f'(t) \{ \sin^2(\theta)  -    \sin^2(\theta_j)\} dt$$
By using lemma 2.1 we have
$$
I_1\leq  C_1\lambda ^{-\frac{1}{2}} \sup_{ a\leq x\leq b}|\theta(x)-\theta_{j}(x)|.
$$
Applying lemma 2.1 and lemma 2.2 we have
\begin{eqnarray*}
I_2 &:=&\lambda ^{-\frac{1}{2}} \int_{a}^{x} f^2(t) \{ \sin(2\theta)  -    \sin(2\theta_j)\} \sin (2\theta) dt\\
&+& \lambda ^{-\frac{1}{2}} \int_{a}^{x} f^2(t)\{ \sin(2\theta)  -    \sin(2\theta_j)\} \sin (2\theta_j) dt\\
&+& \lambda ^{-\frac{1}{2}} \int_{a}^{x} f^2(t) \{ \sin(2\theta_j)  -    \sin(2\theta_{j-1})\} \sin (2\theta_j) dt\\
&\leq&  C_2\lambda ^{-\frac{1}{2}} \sup_{ a\leq x\leq b}|\theta(x)-\theta_{j}(x)|
\end{eqnarray*}
 Finally, using lemma 2.1, we conclude
\begin{eqnarray*}
I_3 &:=&\lambda ^{-1} \int_{a}^{x} fF \{ \sin(2\theta)  -    \sin(2\theta_j)\}\sin^2 (\theta) dt\\
&+& \lambda ^{-1} \int_{a}^{x} fF ( \sin(\theta)  -    \sin(\theta_{j-1})) ( \sin(\theta)  +   \sin(\theta_{j-1})) \sin (2\theta_j) dt\\
&\leq&  C_3\lambda ^{-\frac{1}{2}} \sup_{ a\leq x\leq b}|\theta(x)-\theta_{j}(x)|
\end{eqnarray*}
This ends the proof of lemma 2.5.
\vskip 0.4 true cm
\begin{thm} Suppose that (2.3) hold for some positive integer $N$, then
$$
\theta(b)- \theta(a)- (b-a)\lambda^\frac{1}{2}=-\int_{a}^{b}f\sin (2\theta_N(x))dx +  (\lambda^\frac{-1}{2})\int_{a}^{b}F\sin^2 (\theta_N)dx + o(\lambda^{\frac{-N}{2}})
$$
as $\lambda \rightarrow \infty$.
\end{thm}

\vskip 0.4 true cm
\pn{\bf Proof.} We integrate (1.5) over $[a,x]$ and obtain
$$
\theta(x)-\theta(a)= \lambda ^{\frac{1}{2}}(x-a)-\int_{a}^{x}f\sin(2\theta(t))dt+ \lambda ^{-\frac{1}{2}}\int_{a}^{x}F\sin^2(\theta(t))dt
$$
In particular
$$
\theta(b)-\theta(a)= \lambda ^{\frac{1}{2}}(b-a)-\int_{a}^{b}f\sin(2\theta(t))dt+ \lambda ^{-\frac{1}{2}}\int_{a}^{b}F\sin^2(\theta(t))dt
$$
and so,
\begin{eqnarray*}
\theta(b)- \theta(a)- (b-a)\lambda^\frac{1}{2}&=& -\int_{a}^{b}f\sin (2\theta_N(x))dx +  (\lambda^\frac{-1}{2})\int_{a}^{b}F\sin^2 (\theta_N)dx\\
&+&\int_{a}^{b}f \{ \sin (2\theta_N(x)-\sin (2\theta(x)) )  \}dx\\
&+&(\lambda^\frac{-1}{2})\int_{a}^{b}F \{ \sin^2 (\theta)-\sin^2 (\theta_N)\}dx.
\end{eqnarray*}
We need to prove that two last terms are $o(\lambda^{\frac{-N}{2}})$ as $\lambda \rightarrow \infty$. Applying lemma 2.2 and lemma 2.4 we have

\begin{eqnarray*}
I &:=&\int_{a}^{b}f(x) \{ \sin (2\theta_N(x)-\sin (2\theta(x)) )  \}dx+ (\lambda^\frac{-1}{2})\int_{a}^{b}F(x) \{ \sin^2 (\theta)-\sin^2 (\theta_N)\}dx\\
&\leq& C\lambda ^{-\frac{1}{2}} \sup_{ a\leq x\leq b}|\theta(x)-\theta_{N}(x)|+ C(\lambda^\frac{-1}{2})\int_{a}^{b}F \sup_{ a\leq x\leq b}|\theta(x)-\theta_{N}(x)| dx
\end{eqnarray*}

When $N=1$, applying lemma 2.5, $|\theta(x)-\theta_{1}(x)|= o(\lambda^{\frac{-1}{2}})$. Now By using lemma 2.3 and induction we achieve that $I =o(\lambda^{\frac{-N}{2}})$ as $\lambda \rightarrow \infty$.

\begin{rem}
By using the discussions of choice of $f$ in [10], the condition (2.3) let us to consider $q$ as the form  $q(x)\sim x^{-K}$ where $ 1 \leq K<2$.
\end{rem}


\vskip 0.4 true cm


\bigskip
\bigskip

{\footnotesize \pn{\bf Mahdi Hormozi}\; \\ {Department of Mathematical Sciences, Division of Mathematics}, {Chalmers University of Technology and University of Gothenburg, Gothenburg 41296, Sweden.} \\
{\tt Email:
hormozi@chalmers.se}\\
\end{document}